\documentclass[12pt]{article}

\usepackage[english]{babel}

\usepackage{amsmath}
\usepackage{amsthm}
\usepackage{amsfonts}
\usepackage{amssymb}

\usepackage{stmaryrd}

\usepackage{eucal}
\usepackage{enumerate}
\usepackage{graphicx}

\usepackage[colorlinks=true]{hyperref}

\theoremstyle{plain}
\newtheorem*{thma}{Theorem A}
\newtheorem*{thmb}{Theorem B}

\newtheorem{thm}{Theorem}[section]
\newtheorem{prop}[thm]{Proposition}
\newtheorem{lem}[thm]{Lemma}
\newtheorem{cor}[thm]{Corollary}
\newtheorem{defn}[thm]{Definition}
\newtheorem{defs}[thm]{Definitions}

\theoremstyle{remark}

\newcommand{\set}[1]{\left\{#1\right\}}
\newcommand{\minus}{\smallsetminus}

\newcommand{\tree}{\CMcal{T}}
\newcommand{\cB}{{\CMcal{B}}}

\newcommand{\bC}{\mathbb{C}}
\newcommand{\bR}{\mathbb{R}}

\newcommand{\bZ}{\mathbb{Z}}

\newcommand{\MRef}[2]{\href{http://www.ams.org/mathscinet-getitem?mr=MR{#1}}{MR~{#2}}}

\DeclareMathOperator{\tran}{\textrm{tran}}
\DeclareMathOperator{\Leb}{\textrm{Leb}}

\newcommand{\child}{\textnormal{\textsf{C}}}
\newcommand{\parent}{\textnormal{\textsf{P}}}

\newcommand{\ap}{a^{\parent}}
\newcommand{\ac}{a^{\child}}

\newcommand{\Kjulia}{\CMcal{K}}
\newcommand{\julia}{\CMcal{J}}
\newcommand{\ray}{\CMcal{R}}

\newcommand{\circled}[1]{{$\mathnormal{\bigcirc}$}\hspace{-8.5pt}{\footnotesize
{\textit{#1}}}\hspace{3pt}}

\newcommand{\harm}{\omega}
\newcommand{\Harm}{\Omega}
\newcommand{\dist}[2]{\text{dist}\,(#1,#2)}

\begin{document}

\title{Brownian motion, random walks on trees, and harmonic measure on
  polynomial Julia sets}

\author{Nathaniel D. Emerson}

\maketitle

\begin{abstract}
  We consider the harmonic measure on a disconnected polynomial Julia set in
terms of Brownian motion. We show that with probability one, the first point in the
Julia set that a Brownian particle hits a single-point component. Associated to the
polynomial is a combinatorial model, the tree with dynamics. We define a random walk on
the tree, which is a combinatorial version of Brownian motion in the plane. This random
walk induces a measure on the tree, which is isomorphic to the harmonic measure on the
Julia set.
\end{abstract}

\section{Introduction}\label{sect: introduction}

Let $f$ be a polynomial of degree $d \geq 2$.   The  set of all points with a bounded
orbit under $f$,
\[
    \Kjulia_f = \set{z \in \bC : \ \sup_n |f^n(z)| < \infty},
\]
is called \emph{filled -in Julia set of $f$}, where $f^n$ denotes the $n^{\textrm{th}}$
iterate of $f$. The \emph{Julia set of $f$} is the topological boundary of $\Kjulia_f $:
$\julia_f =
\partial \Kjulia_f$. We will only consider $f$ with $\Kjulia_f$
disconnected. By a classical result of Fatou and Julia, this means
at least one critical point of $f$ has unbounded orbit.

Imagine a particle moving on the Riemann sphere according to the laws of Brownian motion
\cite{Doob}. Suppose the particle starts at the point at infinity and let $Z(t)$ denote
the position of the particle at time $t$. We call $Z(t)$ a \emph{Brownian path}.
Following A.~Lopes \cite{Lopes} and S.~Lalley \cite{Lalley}, we consider the interaction
of a Brownian particle with a polynomial Julia set. We say an event \emph{almost surely
(a.s.)} occurs if the probability of it occurring is 1.  A Brownian path almost surely
enters $\Kjulia_f$ in a finite amount of time.  Let $t_0 = \inf \set{t: \ Z(t) \in
\Kjulia_f}$. We call $t_0$ the \emph{first entry time} of $Z$ into $\Kjulia_f$, and
$Z(t_0)$ the \emph{first entry point}.

Brownian motion induces a measure on $\Kjulia_f$, which we call $\harm_f$, the
\emph{harmonic measure} of $f$ \cite{Doob}. For $X$ a measurable subset of $\Kjulia_f$,
$\harm_{f}(X)$ is the probability that the first entry point of a Brownian path lies in
$X$. The harmonic measure of $f$, is a Borel probability measure defined on $\Kjulia_f$.
It is $f$-invariant, $f$-ergodic and strongly mixing. Roughly, the harmonic measure
describes the one-dimensional structure of $\Kjulia_f$.  A useful feature of the
harmonic measure is that it has a variety of equivalent definitions, see Theorem
\ref{thm: harmonic mes defs}. We will consider harmonic measure in terms of Brownian
particles for intuitive purposes \cite{Lalley}. For technical purposes, we will define
it in terms of landing external rays \cite{Arsove}.

We consider the interaction of disconnected polynomial Julia sets and Brownian
particles. We consider disconnected polynomial Julia sets.  In particular, we consider
Julia sets that have connected components which are not points, which we refer to as
\emph{island components}.  We call a single-point component a \emph{singleton}. Note the
island components are clearly visible in Figure \ref{fig: disc Julia}.  A disconnected
polynomial Julia set will always have uncountably many singletons.  So, the the
structure of the Julia sets that we consider is something like an island chain
surrounded by a barrier reef.  The island components are larger than singletons in a
topological sense.  For instance they have positive diameter.  They are also larger in
terms of Brownian particles, a Brownian path will almost surely visit a given island
component.

\begin{thma}
  If $f$ is a polynomial with a disconnected Julia set, then the first entry point of a
  Brownian path into $\Kjulia_f$ is almost surely a singleton.
\end{thma}

So in terms of Brownian particles, island components are no larger than singletons. This
result implies a number of facts about the structure of Julia sets that have countably
many connected components that are not points.

We give an estimate on how quickly the harmonic measure decreases with respect to
equipotentials of Green's function (Theorem \ref{thm: decay}).  We show that the
harmonic measure always decreases exponentially. It follows that the harmonic measure of
any component of a disconnected polynomial Julia set is zero.

The major tool in this paper is the combinatorial system of a \emph{tree with dynamics}
\cite{E03}. Associated to a polynomial with a disconnected Julia set is a canonical tree
with dynamics. The tree with dynamics is a discrete model for the dynamics of such a
polynomial. It captures many important facets of the dynamics of a polynomial, but is
easier to work with than the polynomial itself. We construct the tree by decomposing the
basin of attraction of the point at infinity into conformal annuli using Green's
function. These annuli have a natural tree structure, which is compatible with the
dynamics.  The polynomial has a well-defined degree on each annulus. We associate each
annulus to a vertex of the tree.  We obtain a countable, rooted tree $\tree$,  a map
$F:\tree \to \tree$, and a degree function $\deg: \tree \to \bZ^+$.  We define \emph{the
combinatorial harmonic measure} on $\tree$.  We put a point-mass at the root of $\tree$.
We distribute the mass of a vertex to its pre-images under $F$, weighting by the degree
of a pre-image.  We show that the harmonic measure of a component of the Julia set can
be estimated by the combinatorial harmonic measure on the tree.  This result allows us
to prove our other main theorem: the combinatorial harmonic measure is a model for the
harmonic measure in the plane.

\begin{thmb}
  If $f$ is a polynomial with a disconnected Julia set, then the combinatorial harmonic measure on the
  tree with dynamics of $f$ is isomorphic to the harmonic measure on the Julia set of $f$.
\end{thmb}

This allows us to use techniques from the filed of discrete
potential theory to answer questions about the harmonic measure.
Harmonic measure in the plane can be defined by of Brownian motion.
The combinatorial harmonic measure can be defined a random walk on
the tree.

The rest of this paper is organized as follows.

In Section \ref{sect: background}, we give some background on
potential theory.  We describe the decomposition of the plane using
Green's function. While this is a standard technique, there are some
subtle points we later use. We recall some facts about the harmonic
measure.

We state our results for Julia sets in Section \ref{sect: Julia Results}. Especially an
estimate on the rate of decrease of harmonic measure. We discuss various consequences of
this result. To show Theorem A, we use a recent result of W.~Qiu and Y.~Yin
\cite{Qiu-Yin}.

We describe the tree with dynamics of a polynomial, in Section
\ref{sect: TwD}. We define the combinatorial harmonic measure on the
tree with dynamics. We give the combinatorial results that imply
Theorem A. We show that the combinatorial harmonic measure is
isomorphic to the harmonic measure on the filled-in Julia set,
proving Theorem B.

While completing this paper, the author learned that L. DeMarco and C. McMullen had
independently obtained many of the same results, in particular Theorems A and B
\cite{DeMarco-McMullen}.


\section{Background} \label{sect: background}

We will consider the interaction of holomorphic dynamics and
potential theory. We assume basic familiarity with holomorphic
dynamics \cite{Steinmetz}. We use two objects from potential theory:
Green's function and harmonic measure. The book of N.~Steinmetz has
a short introduction to Green's function and harmonic measure on
polynomial Julia sets \cite{Steinmetz}. T.~Ransford has written a
very readable introduction to potential theory, which includes a
section on polynomial dynamics \cite{Ransford}. The work of J.~Doob
gives a more complete account of potential theory, and covers
Brownian motion in detail \cite{Doob}.  The paper of S.~Lalley is a
good introduction to the particulars of Brownian motion and Julia
sets of rational functions \cite{Lalley}.

We give the details of the dynamical decomposition of the plane. Following Branner and
Hubbard, we use equipotentials of Green's function of a polynomial to decompose the
plane into conformal annuli \cite{Branner}.  These annuli have a natural tree structure,
which is compatible with the dynamics.  So we can associate a polynomial to the
combinatorial system of tree with dynamics due to R. P{\'e}rez-Marco \cite{E03}.

We then recall some facts about the harmonic measure of a polynomial
Julia set. We consider harmonic measure in terms of landing external
rays \cite{Arsove} and Brownian motion \cite{Lalley}.  We give a
variety of equivalent definitions of harmonic measure.  We discuss
subsets that are \emph{shielded} from the harmonic measure---a
phenomenon that we later show occurs in polynomial Julia sets.

\subsection{The Dynamical Decomposition of the Plane}

We define an \emph{annulus} as a subset of the complex plane that is conformally
equivalent to a set of the form $\set{z \in \bC: \ r_1 < |z| < r_2}$, for some $r_1,
r_2$ with $0 \leq r_1< r_2 \leq \infty$. We say a set $S$ is \emph{nested} inside an
annulus $A $, if $S$ is contained in the bounded components of $\bC \minus A$. For an
annulus $A $, we define the \emph{filled-in annulus}
\[
    P(A) = A \cup \set{\text{bounded components of } \bC \minus A}.
\]
Observe that $P(A)$ is an open topological disk.

For the remainder of this paper, let $f$ be a polynomial of degree
$d \geq 2$ with disconnected Julia set. Let $g$ denote Green's
Function of $f$. The functional equation $g(f) = d \cdot g$ is
satisfied by $f$ and $g$. We use $g$ to define the dynamic
decomposition of the basin of attraction of infinity for $f$.

An \emph{equipotential} is a level set of $g$; $\set{z \in \bC: \
g(z) = \lambda >0}$. The critical points of $g$ are the critical
points of $f$ and the pre-images of critical points of $f$. We
distinguish all equipotentials of $g$ that contain a critical point
of $g$ or an image under $f$ of a critical point of $g$. There are
countably many such equipotentials, say $\set{E_l}_{l \in \bZ}$. We
index them so that $g|E_l < g|E_{l-1}$, $E_l$ is a Jordan curve for
$l \leq 0 $, and $E_1$ is not a Jordan curve (so it contains a
subset homeomorphic to a figure-8). Let $H$ be the number of orbits
of $\set{E_l}_{l \in \bZ}$ under $f$. If  $f$ has $e$ distinct
critical points that escape to infinity, then $H \leq e$. It is
possible that $H<e$, if $f$ has two critical points $c_1$ and $c_2$
such that $g(c_1) = d^n g(c_2)$ for some $n \in \bZ$. From the
functional equation and the indexing of $E_l$, it follows that
$f(E_l) = E_{l-H}$ for all $l$.

Define $ U_l = \set{z: \ g|E_{l} > g(z) > g|E_{l+1}}$. For $l \leq 0$, $U_l$ is a single
annulus.  For all $l$, $U_l$ is open and consists of the disjoint union of finitely many
annuli $A_{l,i}$. We call each of the $A_{l,i}$ an \emph{annulus of $f$} at level $l$.
The closure of a filled-in annulus, $\overline{P}(A_{l,i})$, is called a \emph{puzzle
piece} of $f$ at depth $l$ \cite{Branner}. A sequence $(A_l)_{l=0}^{\infty}$ of annuli
of $f$ is called \emph{nested}, if $A_l \subset U_l$ and $A_{l+1}$ is nested inside $A_l
$ for all $l$. The intersection of the filled-in annuli $\bigcap_{l=0}^{\infty} P(A_l)$
from nested sequence is a connected component of $\Kjulia_f$.

We code the dynamics of a polynomial with a disconnected Julia set by the combinatorial
system of a \emph{tree with dynamics} \cite{E03}.  We use $\set{A_{l,i}}$, the annuli of
$f$, to form a tree $\tree$ by associating  each $A_{l,i}$ to a vertex $a_{l,i}$. We
define a map $\tau: \tree \to \set{A_{l,i}}$ by $\tau(a_{l,i}) = A_{l,i}$.  We declare
that there is an edge between $a_{l,i}$ and $a_{l-1,j}$ if $A_{l,i}$ is nested inside
$A_{l-1,j}$.  From the functional equation $g(f) = d \cdot g$, we can show that the
image of any annulus of $f$ is another annulus of $f$. That is, for any $A_{l,i}$, we
have $f(A_{l,i}) = A_{l-H, j}$ for some $j$. So the dynamics are compatible with the
tree structure, and we define $F: \tree \to \tree $ by
\[
    F(a) = b \quad \text{if } \quad
    f(\tau(a)) = \tau(b).
\]
Note that $\tau$ conjugates $F$ to $f$, that is $\tau(F(a)) = f(\tau(a))$. We define the
\emph{degree} of each vertex $\deg a_{l,i}$ as the topological degree of $f|A_{l,i}$. We
call the triple $<\tree, F, \deg>$ the \emph{tree with dynamics} of $f$. When drawing
trees, we show a vertex of degree 1 by $\bullet$, and a vertex of degree $ D > 1$ by
\circled{D}.

\begin{figure}[hbt] \label{fig: cubic R=1}
\begin{center}
\includegraphics{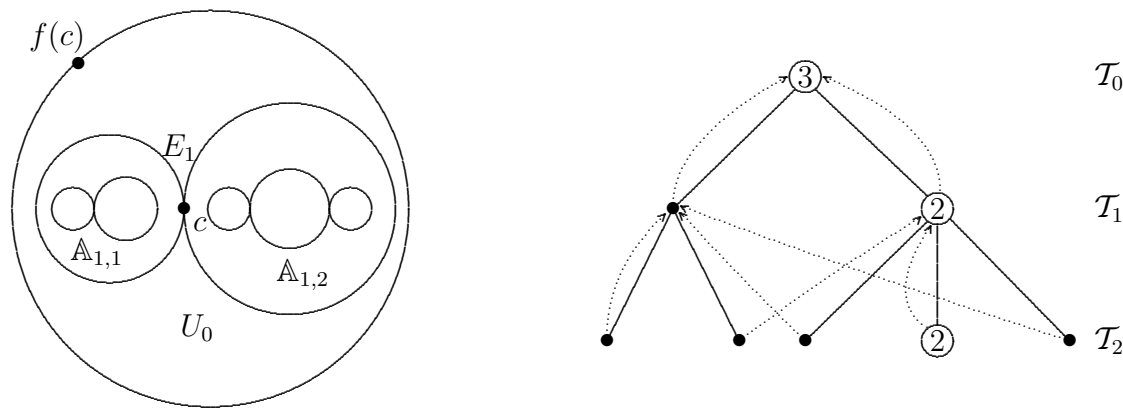}
\caption{Equipotentials of a cubic polynomial with one escaping
critical point, $c$ (left).  Its tree with dynamics, with $H=1$
(right).}
\end{center}
\end{figure}

\subsection{Harmonic Measure}

There is a Borel probability measure on $\Kjulia_f$, the
\emph{harmonic measure} $\harm_f$ \cite{Ransford}. In fact, the
support of $\harm_f$ is always contained in $\julia_f$. So whether
one considers $\harm_f$ a measure on $\Kjulia_f$ or $\julia_f$ is a
matter of preference. The harmonic measure is $f$-invariant,
ergodic, and non-atomic. The harmonic measure is always mutually
singular to two-dimensional Lebesgue measure \cite{Oksendal}. The
support of $\harm_{f}$ has Hausdorff dimension at most 1
\cite{Jones-Wolff}.  The harmonic measure depends only on the
topology of $\Kjulia_f$, and not the conformal structure.

For technical purposes we will define the harmonic measure in terms
of landing external rays.  This is a special case of the Green's
measure \cite{Brelot-Choquet}. For intuitive purposes, we will
consider the harmonic measure in terms of Brownian motion
\cite{Lalley}.

A \emph{Green's line} is a is an orthogonal trajectory to the equipotentials of Green's
function.  For a polynomial of degree at least 2, each Green's line can be canonically
identified with a point in the circle at infinity, that is the set of asymptotic
directions in the plane. An \emph{external ray}, $\ray_{\theta}$, is a Green's line
labelled with an angle $\theta$ from the circle at infinity $\mathbb{T}  = \bR / \bZ $.
We say an external ray is \emph{smooth} if it does not contain a critical point of
Green's function.   All but countably many external rays are smooth. Hence in terms of
Lebesgue measure on the circle, almost all external rays are smooth. A smooth external
ray intersects each equipotential of Green's function in a unique point, so we can
parameterize an smooth ray by potential. For $\lambda >0$, let $\ray_{\theta}(\lambda)$
be the the unique point in $\ray_{\theta} \cap \set{g = \lambda }$.  We say that a
smooth ray $\ray_{\theta}$ \emph{lands} at $z \in \Kjulia_f$ if $\lim_{\lambda \to 0^+}
\ray_{\theta}(\lambda) = z$.  Almost every external ray lands. For $X$ a measurable
subset of $\Kjulia_f$, we have $\harm_{f}(X) = \Leb_1 (\set{\theta: \ \ray_{\theta}
\text{ lands at } X})$ \cite{Brelot-Choquet, Arsove}, where $\Leb_1 $ denotes the
normalized Lebesgue measure on the unit circle $\mathbb{T}$.

There are a variety of characterizations of the harmonic measure
which we recall now.

\begin{thm} \label{thm: harmonic mes defs}

 If $\Kjulia$ is the filled-in Julia set of a polynomial of degree at least 2, then the
following measures on $\Kjulia$ are equal:
\begin{enumerate}[{\indent}a.]
    \item the harmonic measure  \cite{Ransford};
    \item the equilibrium measure  \cite{Ransford};
    \item the Green's measure  \cite{Arsove};
    \item the hitting measure of Brownian motion \cite{Doob};
    \item the Brolin measure  \cite{Brolin}.
\end{enumerate}
In fact, for $\Kjulia$ a compact subset of $\bC$ with positive capacity, a--d are always
equal.
\end{thm}

\begin{center}
\begin{figure}[hbt] \label{fig: hidden harmonic}
\includegraphics{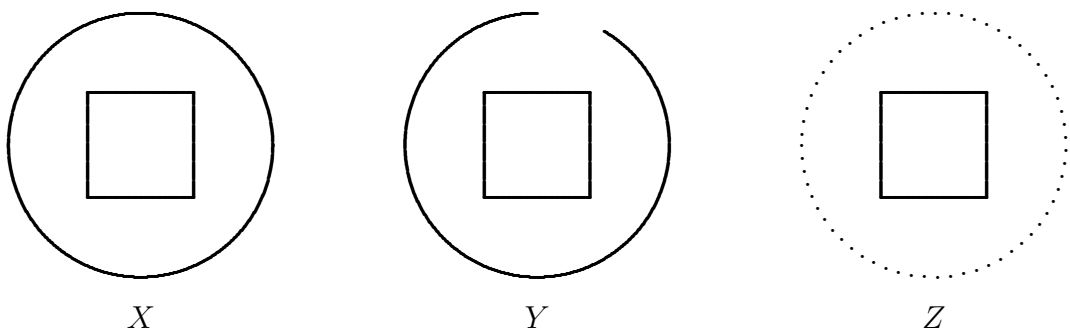}
\caption{A subset shielded from the harmonic measure.}
\end{figure}
\end{center}

We briefly consider the harmonic measure of more general compact
subsets of the plane. If $X \subset \bC$ is compact with positive
capacity, the harmonic measure $\harm_X$ (with a pole at the point
at infinity) is defined on $X$ \cite{Ransford}.  If $X$ is a
rectifiable curve, then $\harm_X$ is just the normalized
one-dimensional Lebesgue measure on $X$. For instance, if $X$ is a
circle and $A$ an arc of angle $\theta$, then $\harm_X (A) =
\theta/2 \pi$. If $X$ is a square and $S$ one of its sides, $\harm_X
(S) = 1/4$. The harmonic measure of $C \subset X$ depends not only
on the intrinsic properties of $C$, but on how $C$ is embedded in
$X$. If $A \cup B$ is a partition of $X$, let us say $A$
\emph{shields} $B$ if the first entry point of a Brownian particle
almost surely lies in $A$. That is, if $\harm_X(A) =1$.   There is
an intuitive explanation for shielding in terms of Brownian motion.
Refer to Figure \ref{fig: hidden harmonic}.  Let $S$ be a square and
$C$ be a circle.  Take $X$ as $S$ enclosed by $C$. For this $X$,
$\harm_X (S) = 0$, since a Brownian path must intersect the circle
before it hits $S$.  Now let $Y$ be $S$ enclosed by $C$ and remove
an arc of angle $ \theta $ from $C$. Then $\harm_Y(S) \leq \theta /
2 \pi$, since a Brownian path whose first entry point is in $S$,
must first pass through the gap in $C$. Now imagine we form $Z$ by
removing countably many arcs from $C$, so that what remains is a
Cantor set of length $L$, where $0 \leq L <1$. It follows that
$\harm_Z(S) \leq 1 - L$.

We show that shielding occurs for Julia sets. If $\Kjulia_f$ is a disconnected
polynomial Julia set with a component that is not a singleton, then the singletons of
$\Kjulia_f$ shield the island components, see Theorem A.  This is similar to the
situation in the set $Z$ above.  In light of Theorem \ref{thm: harmonic mes defs}, the
islands components are also shielded from external rays.


\section{Results for Julia sets} \label{sect: Julia Results}

In this section, we give results for disconnected polynomial Julia
sets. We only consider disconnected Julia sets.  We will assume two
results from Section 4, Theorem \ref{thm: tree decay} and Lemma
\ref{lem: Harm(a) = harm(A)}.

First, we will consider arbitrary disconnected Julia sets.  We
extend the harmonic measure to annuli of $f$.  We show that the
harmonic measure of an annulus decreases as an exponentially with
the level of the annulus.

We then restrict our attention to disconnected Julia sets with a component that is not a
point. H.~Brolin gave the first example of such a Julia set \cite[p.137--138]{Brolin}.
B.~Branner and J.~Hubbard showed if the component containing a critical point of $f$ is
periodic, then that component and all of its pre-images are not singletons \cite[Thm.
5.3]{BH92}. Recently, Qiu and Yin announced in a pre-print that $\Kjulia_f$ is a Cantor
set unless it has a periodic component containing a critical point \cite[Main
Thm.]{Qiu-Yin}.

\subsection{Rate of decrease for the harmonic measure}

We now extend the notion of harmonic measure to an annulus of $f$.  By slight abuse of
notation, we use $\harm_f$ to represent this notion.

\begin{defn}
Let $A$ be an annulus of $f$.  Define
\[
    \harm_f(A) = \Leb_1 \set{\theta : \ \CMcal{R}_{\theta} \cap
    A \neq \emptyset}.
\]

\end{defn}

We give an estimate on the harmonic measure of components of the Julia set in terms of
$\harm_f(A)$.

\begin{lem}\label{lem: approx omega}
Let $f$ be a polynomial of degree $\geq 2$. Let $K$ be a component of $\Kjulia_f$.  Let
$\set{A_l}_{l=0}^{\infty}$ be the unique nested sequence of annuli of $f$ such that $K =
\bigcap_{l=0}^{\infty} P(A_l)$. Then
\[
    \harm_{f}(K) \leq \lim_{l \to \infty} \harm_f(A_l).
\]
\begin{proof}
By Theorem \ref{thm: harmonic mes defs}, $\harm_f(K)$ is the measure of external rays
that land on $K$.  For each $l$, $K$ is nested inside $A_l$  and the boundary of $A_l$
is contained in two equipotentials. Hence, any ray that lands on $K$ must also intersect
$A_l$. Thus, $\harm_{f}(K) \leq \harm_f(A_l)$ for each $l$. Since $A_{l+1} $ is nested
inside $A_l$, we have $ \harm_f(A_{l+1}) \leq \harm_f(A_l)$. Taking a limit finishes the
lemma.

\end{proof}
\end{lem}

Lemma \ref{lem: approx omega} is an important result for proving Theorem A.  It allows
us to estimate $\harm_f(A_l)$, instead of computing $\harm_f(K)$ directly.  We can
easily transfer these estimates to the tree with dynamics.  We show $\harm_f(A_l)$
decreases exponentially with $l$.  This follows from the analogous combinatorial result,
Theorem \ref{thm: tree decay} and Lemma \ref{lem: Harm(a) = harm(A)}. Let $\lceil \cdot
\rceil$ denote the ceiling function.

\begin{thm}\label{thm: decay}
Let $f$ be a polynomial of degree $d \geq 2$ with a disconnected Julia set. Let $D= 1 +
M$, where $M$ is the maximum of the multiplicities of the non-escaping critical points
of $f$.   Let $H$ be the number of orbits of escaping critical points under $f$. There
exists a constant $c_0 > 0 $, such that if $A$ is an annulus of $f$ at level $l$, then
\[
    \harm_f(A)  \leq
    c_0\left( \frac{D}{d} \right)^{\lceil l/H \rceil}.
\]

\end{thm}

We restate the above theorem in terms of escaping critical points.
This is of interest if one considers a polynomial where the number
or multiplicity of escaping critical points is known, but $H$ is
not. For instance, in the case of a polynomial from some escape
locus in parameter space \cite{BH92}.

\begin{cor}
Let $f$ be a polynomial of degree $d \geq 2$.  Let $D$, $H$, and $A$ be the same as in
Theorem \ref{thm: decay}.  Let $e$ be the number of distinct critical points of $f$ that
escape to infinity. Let $m$ be the number of critical points of $f$, counted by
multiplicity, that escape to infinity. Then

\[
    \harm_f(A)
    \leq c_0\left( \frac{D}{d} \right)^{\lceil l/H \rceil}
    \leq c_0 \left( \frac{D}{d} \right)^{\lceil  l/ e \rceil}
    \leq c_0 \left( \frac{D}{d} \right)^{\lceil  l/ m \rceil}.
\]
\begin{proof}
We have $H \leq e \leq m$, so the last two inequalities are easily
verified.
\end{proof}
\end{cor}

It follows that the harmonic measure of a component of a disconnected Julia set is 0.

\begin{cor}\label{cor: mes(comp) = 0}
\begin{proof}
For every $l \geq 0$, there is a unique annulus $A_l$ of $f$ at level $l$ such that $K $
is nested inside $A_l$.  Combining Lemma \ref{lem: approx omega} and Theorem \ref{thm:
decay}, we obtain
\[
    \harm_{f}(K) \leq \lim_{l \to \infty} \harm_f(A_l) \leq \lim_{l \to \infty} c_0\left( \frac{D}{d} \right)^{\lceil l/H \rceil}.
\]
Since $\Kjulia_f$ is disconnected, we have $D < d$, so the right hand side tends to 0 as
$l$ approaches $\infty$ .
\end{proof}
\end{cor}

That is to say, no component of $\Kjulia_f$ is charged by $\harm_{f}$.  Note that one
could prove this by $f$-invariance of $\harm_f$.

\subsection{Julia sets with island components}

For the rest of the section we assume that the Julia set has a component that is not a
singleton. If $K$ is a
 component of $\Kjulia_f$ that is not a singleton, it will have positive capacity. Nonetheless, it
will not be charged by $\harm_{f}$.

\begin{center}
\begin{figure}[hbt]\label{fig: disc Julia}
\includegraphics[scale=.4]{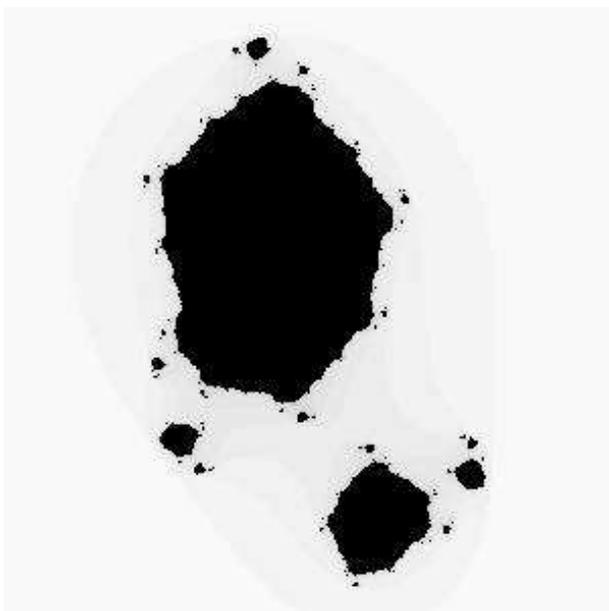}
\caption{A disconnected polynomial Julia set with island
components.}
\end{figure}
\end{center}

We partition the Julia set into singletons and non-singletons.   Let $K(z) $ denote the
connected component of a point $z$ in $\Kjulia_f$. Define
\[
    \Kjulia_f^0 = \set{z: K(z) = \set{z}}
    \quad \text{and} \quad \Kjulia_f^1 = \set{z: K(z) \neq \set{z}}.
\]
In Figure \ref{fig: disc Julia}, the large ``islands'' are the
components of $\Kjulia_f^1$. The points can be thought of as
$\Kjulia_f^0$.  We study the harmonic measure and this partition.
Since the partition is $f$-invariant and $\harm_f$ is ergodic, one
of these sets must have harmonic measure zero.

Any component of $\Kjulia_f^1$ has positive diameter. Thus, in a
topological sense, components of $\Kjulia_f^1$ are larger than
components of $\Kjulia_f^0$. We ask which of $\Kjulia_f^0$ and
$\Kjulia_f^1$ is larger in the sense of harmonic measure? This is
not just a question of intrinsic properties of the components, but
depends on how $\Kjulia_f^0$ and $\Kjulia_f^1$ are embedded in
$\Kjulia_f$.

In terms of  Brownian motion, a natural way to partition the Julia set is into those
components that a Brownian path almost surely does not enter, and those components that
it almost surely does enter. This is the same partition of the Julia set given above:
$\Kjulia_f^0$ is the points of $\Kjulia_f$ that lie in components that the path almost
surely does not enter, and $\Kjulia_f^1$ is the points of $\Kjulia_f$ that lie in
components that the path almost surely enters:
\[
    \Kjulia_f^0 = \set{z: Z(t) \text{ a.s. does not enter } \Kjulia_f(z)}
    \quad \text{and} \quad \Kjulia_f^1 = \set{z: Z(t) \text{ a.s. enters } \Kjulia_f(z)}.
\]
Hence, the components of $\Kjulia_f^1$ are ``larger'' than the components of
$\Kjulia_f^0$ for Brownian motion in some sense. However, we prove the probability that
the first entry point of a Brownian path lies in any given component is zero.  So for
Brownian motion, one could say that no component of $\Kjulia_f$ is larger than any
other. The explanation for this apparently contradictory result is that $\Kjulia_f^0$
shields $\Kjulia_f^1$ from Brownian particles.  We restate Theorem A, in terms of our
partition.

\begin{thma}
  If $\Kjulia_f$ is a disconnected polynomial Julia set, then the first entry point of a
  Brownian path almost surely lies in $\Kjulia_f^0$.

\begin{proof}
By Corollary \ref{cor: mes(comp) = 0}, the measure of every component of $\Kjulia^1_f$
is zero. Qiu and Yin announced that all but at most countably many components of
$\Kjulia_f$ are singletons \cite[p. 26]{Qiu-Yin}. It follows that the harmonic measure
of the non-singletons of $\Kjulia_f$ is zero.
\end{proof}
\end{thma}

That is, a Brownian path almost surely hits a point of $\Kjulia_f^0$ before it hits a
point of $\Kjulia_f^1$.  Although individually the points of $\Kjulia_f^0$ are
negligible, taken together they form an impenetrable barrier.  Intuitively, the
singleton components are rocks that prevent ships from landing on the islands
components.

Theorem A also tells us where external rays land in $\Kjulia_f$. Intuitively, we can say
that $\Kjulia_f^0$ shields $\Kjulia_f^1$ from external rays.

\begin{cor}\label{cor: ae ext ray lands on Cantor}
If $\Kjulia_f$ is a disconnected polynomial Julia set, then almost every external ray
lands on a singleton of $\Kjulia_f$.
\end{cor}

 G. Levin and F. Przytycki have shown that for $\Kjulia_f $
disconnected, if $K$ is a periodic or pre-periodic component of $\Kjulia_f$, then some
external ray land on $K$ and every accessible point $z \in K$ is accessible via an
external ray \cite[]{Levin-Przytycki}. Corollary \ref{cor: ae ext ray lands on Cantor}
can be thought of as a measure theoretic dual of their result. Topologically components
of $\Kjulia_f$ are visible, however they are shielded from the harmonic measure. e

\section{The Tree with Dynamics} \label{sect: TwD}

This section is the technical heart of this paper.  We work with the
tree with dynamics. First, we recall some properties of the tree
with dynamics.  We then define a version of harmonic measure on the
tree and show that it agrees with the harmonic measure on the annuli
of $f$.  We extend the measure to the boundary of the tree.  We show
that the measure on the boundary of the tree is isomorphic to the
harmonic measure in the plane.  We use the measure on the tree to
define a random walk.  We note the equivalence between random walks
on the tree and Brownian motion in the plane

\subsection{Preliminaries}

 We recall some notation from Section \ref{sect: background}.  We decompose the
basin of attraction of infinity into open sets $\set{U_l}_{l \in \bZ}$, bounded by
equipotentials of Green's function.  For each $l$, $U_l = \bigcup_{i=1}^{n}A_{l,i}$,
where each $A_{l,i}$ is an annulus of $f$. We form the tree with dynamics by associating
a vertex of $\tree$ to an annulus of $f$. The map $\tau: \tree \to \set{A_{l,i}}$
witnesses this association.  The tree with dynamics is a triple $<\tree, F, \deg>$,
where $\tree$ is a tree, $F: \tree \to \tree$ is the dynamics, and $\deg : \tree \to
\bZ^+$ is a degree.

We briefly state some properties of a tree with dynamics without proof. A more complete
discussion can be found in a previous paper of the author \cite{E03}.  The tree $\tree$
is a countable.  That is, it is a countable graph with no non-trivial circuits.  It can
naturally be decomposed into levels.

\begin{defs}\label{defn: tree}
For $l \in \bZ$, define $\tree_l = \tau^{-1} (U_l)$.  For $l \leq 0$, $\tree_l$ consists
of a single vertex, say $\tree_l = \set{a_l}$. We call $a_0$ the \emph{root} of $\tree$,
and $\set{a_{-l}}_{l= 1}^{\infty}$ the \emph{extended root} of $\tree$. Let $a \in
\tree_l$ for some $l \in \bZ$.  We call the unique vertex in $\tree_{l-1}$ that is
adjacent to $a$ the \emph{parent of} $a$, and denote it by $\ap $. Any vertex in
$\tree_{l+1}$ that is adjacent to $a$ is called a \emph{child of} $a$, and denoted by
$\ac$.

\end{defs}

Our convention in drawing trees is that a parent is above its children, as in a
genealogic tree.  So $\ap$ is above $a$ and any $\ac$ is below $a$.  We generally denote
the set of all children of $a$ by $\set{\ac}$. When it is necessary to distinguish among
the children of $a$ we use the notation $\set{a^{\child_i}}$.  The structure of the
extended root is trivial.  Its main purpose is to ensure that all iterates of the
dynamics are defined.

\begin{lem} \label{lem: child axioms}
The tree satisfies the following properties:
\begin{enumerate}[{\indent}a.]
  \item For any $a \in \tree$, there is at least 1 child of $a$.  That is, $\tree$ has no
  leaves.
    \item For any $a \in \tree$, there are only finitely many children of $a$.  That is,
  $\tree$ is locally finite.
  \item The root of $\tree$, $  a_0$, has at least 2 children.
\end{enumerate}
\end{lem}

\begin{lem} \label{lem: dyn axioms}
The dynamics satisfy the following properties:
\begin{enumerate}[{\indent}a.]
  \item The dynamics are \emph{children preserving}. For any $a\in \tree$, the image of a child of $a$ is a child
of $F(a) $. Symbolically, $F(\ac) = F(a)^{\child}$.

    \item There exists $H \in \bZ^+$ such that if $a \in
    \tree_l$, then $F(a) \in \tree_{l-H}$.

  \item \label{sublem: local cover property} The dynamics are locally a branched cover of $\tree$.
For any $a\in \tree$, for each child $F(a)^{\mathsf{C}}$ of $F(a)$
we have
  \[
     \sum_{\set{F(a^{\child_i}) = F(a)^{\child}}}
    \deg a^{\child_i} = \deg{a},
    \]
    We refer to this as the \emph{local cover property.}

\end{enumerate}

\end{lem}

\begin{lem}\label{lem: deg axioms}
The degree function satisfies the following properties:
\begin{enumerate}[{\indent}a.]
  \item Then the degree is \emph{monotone}; for all $ a, \in \tree$, if
$a^{\child}$ is a child of $a$, then $\deg a^{\child} \leq \deg a$.
  \item We have $ \deg a_0 =  \deg
    a_{-l} $, for all $l \geq 1$.
    \item We have $\deg a_0 > \deg a$, for all $a \in \tree_l$ with
    $l \geq 1$.

\end{enumerate}
\end{lem}

Throughout this paper, let $\deg a_0 = d$.  We say that $\tree$ is a
tree with dynamics \emph{of degree} $d$.

\begin{center}
\begin{figure}[hbt] \label{fig: Fib tree}
\includegraphics{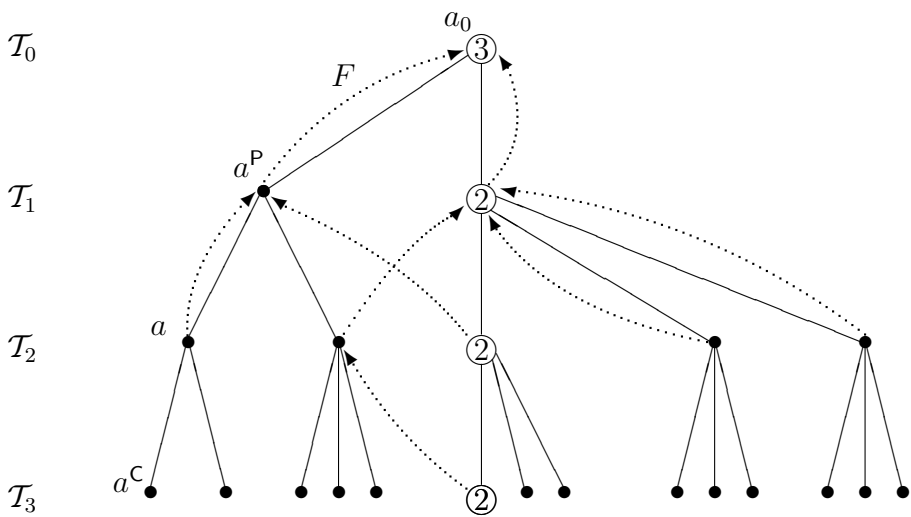}
\caption{A tree $\tree$ with dynamics $F$ of degree 3.  For clarity,
most of the dynamics from level 3 and the extended root are not
shown.}
\end{figure}.
\end{center}

The dynamics are a $d$-fold branched cover of $\tree$ by itself.

\begin{lem}{\cite[Lem. 4.11]{E03}}\label{lem: d pre-images for vertices}
Any vertex of $\tree$ has exactly $d$ pre-images under $F$, counted
by degree. That is, for any $a \in \tree$,
\[
    \sum_{\set{b \in F^{-1}(a) }} \deg b = d.
\]
\end{lem}

We consider all infinite geodesics from the root that move down the
tree.

\begin{defs}
An \emph{end} of ${\tree}$ is a sequence
$\vec{x}=(x_l)_{l=0}^{\infty}$, where $x_l \in \tree_l$ and
$x_{l+1}$ is a child of $x_l$ for all $l$.  Note that a children
preserving map takes ends to ends, so $F(\vec{x})$ is well defined.
We define the degree of an end $\vec{x} = (x_l)_{l=0}^{\infty}$, by
$\deg \vec{x} = \lim_{l \to \infty} \deg x_l$. If $\deg {\vec{x}}  >
1$, then $\vec{x}$ is called a \emph{critical end}.  Let $\cB$
denote the set of all ends of $\tree$.  We call $\cB$ the \emph{end
space of $\tree$}.
\end{defs}

We can define an ultra-metric on $\cB$ by
\[
    \dist{\vec{x}}{\vec{y}} = e^{-L}, \quad \text{where } L = \sup \set{l: x_l =
    y_l},
\]
for ${\vec{x}} \neq {\vec{y}}$, and $\dist{\vec{x}}{\vec{x}} =0$.
 With this metric, $\cB$ is a Cantor set.  Since $F$ is
children-preserving, it extends to a continuous map on $\cB$.  This metric restricts to
a metric on $\tree$, and the end space can naturally be regarded as the boundary of
$\tree$ in this topology.

We can extend $\tau$ in a natural way to a  map from $\cB$ to $\Kjulia_f/\sim$, where
$z_1 \sim z_2$ if $z_1$ and $z_2$ are in the same connected component of $\Kjulia_f$.
Recall, that if $A$ is an annulus, then $P(A)$ is the filled-in annulus: $A \cup
\set{\text{bounded components of } \bC \minus A}$.

\begin{defn}
Let $\vec{x} = (x_l )_{l=0}^{\infty} \in \cB$.  For each $l$, $\tau(x_l)$ is an annulus
of $f$.  Define $\tau: \cB \to \Kjulia_f / \sim$  by
\[
    \tau(\vec{x}) = \bigcap_{l=0}^{\infty} P(\tau(x_l)).
\]

\end{defn}

\begin{prop} \label{prop: eta a homeo}
The map $\tau: \cB \to \Kjulia_f/ \sim$ is a homeomorphism.
\begin{proof}

A Cantor set is homeomorphic to an inverse limit system given by a sequence of
non-trivial open/closed partitions of itself, where the partition at stage $l+1$ is a
refinement of the partition at stage $l$ \cite[Thm. 2--96]{Hocking-Young}.  For $\cB$
one such inverse limit system is given by
\begin{align*}
    \tree_l &\leftarrow \tree_{l+1} \\
    \ap &\mapsfrom a.
\end{align*}

For $\Kjulia/\sim$ the equivalent inverse limit system is given by
\begin{align*}
    U_l &\leftarrow U_{l+1}\\
    A^{\parent} &\mapsfrom A,
\end{align*}
where $A^{\parent}$ is the unique annulus of $f$ at level $l$ that $A$ is nested inside.

By definition, $\tau$ induces an isomorphism of these inverse limit systems.  Therefore
it is a homeomorphism.

\end{proof}
\end{prop}

\subsection{Combinatorial Harmonic Measure}

\begin{defn}
Let $\tree$ be a tree with a distinguished root.  A \emph{flow} on $\tree$ is a function
$\Omega: \tree \to [0, \infty)$ such that for all $a \in \tree$ we have
\[
    \Omega(a) = \sum_{\set{\ac}} \Omega(\ac).
\]

\end{defn}

It is well known that flows on $\tree$ are in one-to-one correspondence with finite
measures on $\cB$.  In Theorem \ref{thm: Extend Harm to ends} we outline the proof of
this fact.

A useful way to think of a flow is in terms of electrical networks \cite{Doyle-Snell}.
Imagine the tree is an electrical network, grounded at its ends, and a charge is
introduced at the root.  The electricity will flow from the root to the ends.  For each
$a \in \tree$, $\Omega(a)$ is current that flows through $a$. Equivalently, one could
imagine that the edge from $\ap$ to $a$ is a wire, and $\Omega(a)$ is its conductance.
The total charge on a set of ends is the measure of the set.

We define $\Harm$, a combinatorial version of harmonic measure. Intuitively, the measure
of $a$ is distributed to the $d$ pre-images of $a$. Each pre-image receives an amount of
measure proportional to its degree.  Compare to the Brolin Measure \cite[\S16]{Brolin}.

\begin{defn} \label{defn: combinatorial measure} For $a \in \tree$, define
$\Harm(a) $ by $\Harm(a) =1$ if $a$ is in the extended root of $\tree$, and
\[
   \Harm(a) = \frac{\deg a}{d}\, \Harm(F(a))
\]
otherwise.  We call $\Harm$ the \emph{combinatorial harmonic
measure} of $\tree$.

\end{defn}

Although we refer to $\Harm$ as a ``measure,'' at this moment $\Harm(a)$ is just a
\emph{weight}---a number associated to each $a \in \tree$. It is not clear that it is a
flow, and since the proof is rather technical we defer it. Nonetheless, we can use
$\Harm$ to estimate $\harm_f(A)$.  We prove Theorem \ref{thm: tree decay}, which implies
Theorem \ref{thm: decay} and thus most of the results in Section \ref{sect: Julia
Results}. We then show that $\Harm$ is a flow. We outline the extension of $\Harm$ to a
measure on $\cB$. Finally, we show that $\tau$ is a measure isomorphism between $(\cB,
\Harm)$ and $(\Kjulia, \harm_f)$.

\begin{center}

\begin{figure}[hbt] \label{fig: Tree w/mes}
\includegraphics{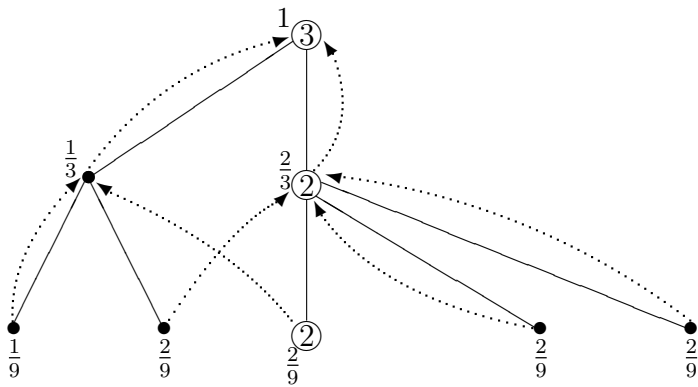}
\caption{A tree  with dynamics with the measure $\Harm (a)$ shown.}
\end{figure}
\end{center}

It is worth noting that by Lemma \ref{lem: d pre-images for
vertices}, $\Harm$ is $F$-invariant.  In the sense that for any $a
\in \tree$,
\[
    \Harm(F^{-1}(a)) = \Harm(a).
\]
Although, we will not use this fact in this paper.

\begin{lem}\label{lem: mes = prod}
If $a \in \tree_l$ for some $l \geq 0$, then
\[
    \Harm (a) = d^{-\lceil l/H \rceil} \prod_{n=0}^{\lceil l/H \rceil} \deg F^n(a).
\]
\begin{proof}
By definition of $\Harm$, we have
\begin{align}
    \Harm(a) &= \frac{\deg a}{d}\, \Harm(F(a))\\
            &= \frac{\deg a}{d}\, \frac{\deg F(a)}{d}\,
            \Harm(F^2(a))\\
            &= d^{-2} \deg a \, \deg F(a) \, \Harm(F^2(a)).
\end{align}
Say that $l = kH +h$, for $0 \leq k$ and $0 < h \leq H$. Then, $
\lceil l/H \rceil = k +1$, so $F^{k+1}(a) = a_{h-H}$.  Repeat the
above argument $k$ times.
\end{proof}

\end{lem}

The measure decreases exponentially with the level of the tree. Note
that this is a combinatorial version of Theorem \ref{thm: decay}.

\begin{thm}\label{thm: tree decay}

Let $D = \max \deg \vec{x}$, for  ends $\vec{x} \in \cB$.  There
exists a constant $c_0$ such that for all $l \geq 0$, if $a \in
\tree_l$, then
\[
    \Harm (a) \leq c_0 \, \left( \frac{D}{d} \right)^{\lceil l/H \rceil}.
\]
\begin{proof}
There are a finite number of levels of the tree that contain a vertex $b$ with $\deg b >
D$, say $Q$ of them, and let $q= \lceil Q/H \rceil$.  Define $c_0 = (d/D)^{q}$. Because
the dynamics go up $H$ levels, the iterates of a point can hit at most $n$ of the levels
with a vertex of high degree.  Hence with at most $q$ exceptions, we have $\deg F^n(a)
\leq D$, so we can replace those terms in Lemma \ref{lem: mes = prod} with $D$. For the
exceptional iterates $\deg F^n(a) \leq d$, and $c_0$ was defined in such a way to
reflect this.

\end{proof}

\end{thm}

The above estimate is sharp when there is an end $\vec{y}$ with $\deg \vec{y} = D$ and
$F(\vec{y}) = \vec{y}$. In general, we can get a little better estimate for a particular
end.

\begin{cor}
Let $\vec{x} \in \cB$, then there exists a constant $c  =
c(\vec{x})$, such that
\[
   \Harm(x_l) \leq c  \, \left( \frac{\deg \vec{x}}{d} \right)^{\lceil l/H \rceil}.
\]
\begin{proof}
Similar to the above lemma.  The key difference is that for $l$
sufficiently large, we have $\deg x_l = \deg \vec{x}$.
\end{proof}

\end{cor}

\begin{cor}
Let $\vec{x} \in \cB$, then
\[
   \lim_{l\to \infty} \Harm (x_l) = 0 .
\]
\begin{proof}
Note that $\deg \vec{x} < d$ by Lemma \ref{lem: deg axioms}.3 and
apply the above corollary.
\end{proof}
\end{cor}

We now establish the first part of the correspondence between the
measures on the tree and in the plane.

\begin{lem} \label{lem: Harm(a) = harm(A)}
For all $a \in \tree$,
\[
    \Harm(a) = \harm_f (\tau(a)).
\]
\begin{proof}
Suppose that $a \in \tree_l$ and use induction on $l$.  If $l \leq 0$, there is only one
annulus at level $l$, so $P(\tau(a)) \cap \Kjulia_f = \Kjulia_f$ and $\Harm(a) = 1=
\harm_f (\tau(a))$.  For $l
> 0$, let $\tau(a) = A$.  We compute $\harm_f(A)$
 in terms of $\harm_f(f(A))$. Now $f$ expands arcs of the circle by a
factor of $d$, so we need to multiply $\harm_f(A)$ by $d$.  However, $f|A$ is a $(\deg
f|A)$-to-one map, so we must divide by $\deg f|A$. Therefore, $\harm_f(f(A))= (d/\deg
f|A) \harm_f (a)$.  Note that $F(a) \in \tree_{l-H}$, so the inductive hypothesis
applies to $F(a)$. It follows that
\begin{align*}
    \Harm(a) &= \frac{\deg a}{d} \, \Harm(F(a))
            &\text{by Definition \ref{defn: combinatorial measure},}\\
            &=  \frac{\deg a}{d} \, \harm_f(\tau(F(a))) &\text{by induction,}\\
            &= \frac{\deg f|A}{d} \, \harm_f(f(\tau(a)))  &\text{by definitions of $\deg a$ and $F$,}\\
            &= \harm_f(\tau(a)). &
\end{align*}

\end{proof}
\end{lem}

Therefore, we can transfer the estimate from Theorem \ref{thm: tree
decay} to a nested sequence of annuli of $f$. Which is exactly the
content of Theorem \ref{thm: decay}.  Therefore, all results in
Section \ref{sect: Julia Results} have now been proven.

\subsection{Extending the Measure to the End Space}

We show that $\Harm$ is a flow.  That is, measure is inherited by
children, as well as pre-images.

\begin{lem} \label{lem: mes A = mes AC}
For all $a \in \tree$,
\[
    \Harm (a) = \sum_{\set{\ac}} \Harm (\ac).
\]
\begin{proof}
Let $a \in \tree_l$.  If $l < 0$, it is clear. For $l=0$, that is $a
= a_0$, we have $\set{a_0^{\child}} = F^{-1}(a_{-H+1})$, since
$\tree_1 = \set{a_0^{\child}} $. By Lemma \ref{lem: d pre-images for
vertices},
\[
    \sum_{\set{a \in F^{-1}(a_{-H+1}) }} \deg a = d.
\]
Recall that $\Harm(a_0) = \Harm(a_{-H+1}) =1$.  We have
\begin{align*}
    \Harm (a_0) &= 1\\
            &= \frac{1}{d} \sum_{\set{a \in F^{-1}(a_{-H+1}) }}
            \deg a\\
            &= \sum_{\set{a_0^{\child}}}
            \frac{\deg a_0^{\child}}{d} \, \Harm (a_{-H+1})\\
             &= \sum_{\set{a_0^{\child}}}
            \frac{\deg a_0^{\child}}{d} \, \Harm (F(a_0^{\child}))\\
            &= \sum_{\set{ a_0^{\child}}} \Harm (a_0^{\child}).
\end{align*}

We use induction on $l > 0$.  Let $a \in \tree_l$ and $F(a) = b$.
Note that $b \in \tree_{l-H}$, so $\Harm (b) = \sum \Harm
(b^{\child_j})$, by the inductive hypothesis.

\begin{align*}
    \Harm (a) &= \Harm (b) \frac{\deg a}{d} \, \\
            &= \frac{1}{d} \sum_{\set{b^{\child_j}}} \Harm (b^{\child_j}) \deg a ,
            \quad \text{by induction}\\
            &= \frac{1}{d} \sum_{\set{b^{\child_j}}}
                \Harm (b^{\child_j})
                \sum_{\set{a^{\child_i} \in F^{-1}(b^{\child_j})}}
                \deg a^{\child_i},\\
                 \intertext{by \ref{lem: deg axioms}.\ref{sublem: local cover property} applied to each child $b^{\child_j}$ of $b$,}
            &= \sum_{\set{a^{\child_i}}}
                    \frac{\Harm(F(a^{\child_i}))}{d} \, \deg a^{\child_i}, \quad \text{since }
                    F(\set{a^{\child_i}})  =  \set{b^{\child_j}} \\
             &= \sum_{\set{a^{\child_i}}}
               \Harm (a^{\child_i}).
\end{align*}
\end{proof}
\end{lem}

\begin{defs}
For $a \in \tree$,  we define  $\CMcal{U}_a$, the \emph{cone of
$a$}, as the set of all ends that pass through $a$. That is,
\[
    \CMcal{U}_a = \set{\vec{x} \in \cB: \ a \in \vec{x}}.
\]
Define the measure of the cone of $a$ by
\[
    \Harm(\CMcal{U}_a) = \Harm(a).
\]
\end{defs}

For any $a \in \tree$, $\CMcal{U}_a$ is an open ball in $\cB$; it is
also compact.  The set of all cones is a sub-basis for the topology
of $\cB$.  Moreover, it is an algebra.

Following Cartier, we outline the proof that $\Harm$ is to a measure on $\cB$.

\begin{thm}\cite[Thm. 2.1]{Cartier} \label{thm: Extend Harm to ends}
We can extend $\Harm$ to a complete Borel measure on $\cB$.  We call $\Harm$ the
\emph{combinatorial harmonic measure on $\cB$}.
\begin{proof}
It follows from  Lemma \ref{lem: mes A = mes AC} that $\Harm$ is finitely additive on
cones. Thus, $\Harm$ is a pre-measure. By standard techniques (Carath\'{e}odory's
Theorem), we can extend $\Harm$ to an outer measure and then a measure on $\cB$.
\end{proof}
\end{thm}

In general, a measure induced by a flow is called a harmonic measure. A tree with
dynamics has a preferred harmonic measure.

We now prove that the harmonic measure on the tree is isomorphic to
the harmonic measure in the plane.

We restate Theorem B in more detail.

\begin{thmb}
Let $f$ be a polynomial with disconnected Julia set.  The harmonic measure $\harm_f$ and
the combinatorial harmonic measure $ \Harm$ are isomorphic. Moreover, $\tau$ induces a
measure isomorphism of $(\cB, \Harm)$ and $(\Kjulia_f, \harm_f)$.
\begin{proof}
Let $z_1 \sim z_2$ if they are in the same component of $\Kjulia_f$. Let $\pi :
\Kjulia_f \to \Kjulia_f/ \sim$ be the projection map. We can consider $\pi^* \harm_f$,
the push-forward by the projection map of $\harm_{f}$.  By Theorem A, $\pi$ is a
bijection, except on a set of measure zero. Thus, for all $X \subset \Kjulia_f$
measurable, we have
\[
    \harm_f(X) = \harm_f(\pi^{-1} \pi (X)).
\]
That is, $\pi $ is a measure isomorphism between $(\Kjulia,
\harm_f)$ and $(\Kjulia/\sim, \pi^* \harm_f)$.

We can also consider $\tau^* \Harm$ the push-forward of $\Harm$ defined on
$\Kjulia/\sim$. We have two measures defined on $\Kjulia_f / \sim$, we show that they
are equal.  By Proposition \ref{prop: eta a homeo}, $\tau$ is a homeomorphism.  By Lemma
\ref{lem: Harm(a) = harm(A)}, $\pi^* \harm_f$ and $\tau^* \Harm$ agree on a sub-basis
for the topology of $\Kjulia_f$.  Therefore, they are equal. That is, if $\CMcal{X}
\subset \cB$ is measurable, then

\[
    \Harm(\CMcal{X}) = \harm_f(\pi^{-1} \tau(\CMcal{X})) .
\]

\end{proof}

\end{thmb}

This isomorphism  gives a new method to compute the harmonic measure
of subsets of $\Kjulia_f$ using annuli of $f$.

\begin{thm}
Let $\Kjulia_f$ be a disconnected polynomial Julia set. Let $X$ be a measurable subset
of $\Kjulia_f$.  For $l \geq 0$, let $\set{A_{l,1}, \dots, A_{l,I(l)}}$ be the annuli of
$f$ at level $l$ such that $X \cap P(A_{l,i}) \neq \emptyset$ for $i =1, \dots, I(l)$.
Then
\[
    \harm_f(X) = \lim_{l \to \infty} \sum_{i=1}^{I(l)}
    \harm_f(A_{l,i}).
\]

\begin{proof}
The inequality
\[
    \harm_f(X) \leq \lim_{l \to \infty} \sum_{i=1}^{I(l)}
    \harm_f(A_{l,i}),
\]
follows easily from Lemma \ref{lem: approx omega}.  The opposite inequality is not
clear. We can consider, $\CMcal{X} = \tau^{-1}(X)$ and show the analogous inequality for
$\Harm$.  Filled-in annuli are analogous to cones.  So what we want to show is
\[
    \Harm(\CMcal{X}) = \lim_{l \to \infty} \sum_{i=1}^{I(l)}
    \Harm(a_{l,i}),
\]
where $\CMcal{X} \cap \CMcal{U}_{a_{l,i}} \neq \emptyset$.

We may assume that $\CMcal{X}$ is compact, since $\Harm$ is Borel.
Fix $\varepsilon > 0 $. We can find $\CMcal{V} \subset \cB $ open,
such that $\CMcal{X} \subset \CMcal{V}$ and $\Harm(\CMcal{V}) \leq
\Harm(\CMcal{X}) + \varepsilon$. Cones are open balls in $\cB$ and
$\CMcal{X}$ is compact, so we can find finitely many cones,
$\CMcal{U}_{b_1}, \dots,\CMcal{U}_{b_J} $, such that
\[
    \CMcal{X} \subset \bigcup \CMcal{U}_{b_j } \subset \CMcal{V}.
\]
By replacing ${b_j}$ by $\set{b_j^{\child}}$, several times if
necessary, we may assume that there is an $L$, such that $b_j \in
\tree_L$ for all $j$.  Hence, $\set{a_{L,i}} \subset \set{b_j}$.
Therefore,
\[
    \sum_{i=1}^{I(L)} \Harm(a_{L,i}) \leq
    \sum_{j=1}^{J} \Harm(b_j) \leq
    \Harm(\CMcal{V}) \leq
    \Harm(\CMcal{X}) + \varepsilon.
\]

\end{proof}

\end{thm}


Just as the harmonic measure in the plane can be defined in terms of Brownian particles,
the combinatorial harmonic measure can be defined in terms of a random walk. A
\emph{random walk} on a tree is a discrete time Markov chain on the tree. We imagine a
particle moving around the tree, with the Markov chain describing its position at each
time.  A random walk is defined its transition function $\tran(x,y)$, which gives the
probability that a particle in $x$ will move to $y$. A \emph{nearest neighbor} random
walk is a random walk where $\tran(x,y) = 0$, unless $x$ is adjacent to $y$. A flow
induces a nearest neighbor random walk. The nearest neighbor random walk induced by
$\Omega$ is the following.

\begin{defn} \label{defn: random walk from Harm}
Define a (nearest neighbor) random walk on $\tree$ by
\[
    \tran(a, \ap) = \frac{1}{2}, \quad \tran(a, \ac) =
    \frac{1}{2}\frac{\Harm(\ac)}{\Harm(a) },
\]
and define all other transition probabilities to be zero.
\end{defn}

It is straightforward to show that the above random walk almost
surely hits the end space. That is, it is \emph{transitive}, it
almost surely visits a given vertex a finite number of times. A
\emph{loop-erased} random walk is a random walk without repeated
vertices. We can transform the above random walk into a loop-erased
random walk. We define
\[
   \tran(a, \ac) =
    \frac{\Omega(\ac)}{\Omega(a)},
\]
and all other transition probabilities are zero. Effectively, this
gives a random end $\vec{x} \in \cB$. The above random walks
correspond to the combinatorial harmonic measure.  Given $\CMcal{X}
\subset \cB$ (measurable), $\Harm(\CMcal{X} )$ is the probability
that a (loop-erased) random walk hits $\CMcal{X} $.

  This can be regarded as a combinatorial version of Brownian motion in the plane.
Suppose a Brownian particle starts in $U_0$ and moves randomly in the plane.  It known
that the hitting probability of this Brownian path is still the harmonic measure
\cite[Prop. 9]{Lalley}. Thus, we almost surely obtain a sequence $(A_i)_{i=0}^{\infty}$
of annuli of $f$, where $A_0 = U_0$ and $A_{i+1}$ is the first annulus of $f$ visited by
the Brownian path after $A_i$.  We call this sequence an \emph{itinerary} of the
Brownian path.  We can also obtain a \emph{loop-erased} itinerary by deleting
repetition.  That is, an itinerary $(A_l)_{l=1}^{\infty}$, where $A_l$ is an annulus at
level $l$.  The probability that a given set of itineraries occurring is equal to the
measure of the analogous set of ends.

\noindent{University of Southern California\\ Los Angeles, CA 90089}\\

\noindent E-mail: \href{mailto:nemerson@usc.edu}{\texttt{nemerson@usc.edu}}


\end{document}